\documentclass[12pt]{article}
\pdfoutput=1

\newlength{\abstractwidth}
\usepackage{epsf}
\usepackage{color}
\usepackage{graphicx}

\usepackage{amsmath}
\usepackage{amssymb}
\usepackage{latexsym}
\usepackage{amssymb,amscd,amsmath,amsthm,url}

\usepackage{tikz, pgf}
\usepackage{tkz-fct}
\usepackage{pgfplots}
\usetikzlibrary{shapes.misc}
\usetikzlibrary{shapes,snakes}
\usetikzlibrary{decorations.pathmorphing}	
\usetikzlibrary{decorations.markings}

\renewcommand{\thefootnote}{\fnsymbol{footnote}}
\renewcommand{\thanks}[1]{\footnote{#1}}
\newcommand{\starttext}{
\setcounter{footnote}{0}
\renewcommand{\thefootnote}{\arabic{footnote}}}

\newcommand{\bea}{\begin{eqnarray}}
\newcommand{\eea}{\end{eqnarray}}
\newcommand{\be}{\begin{eqnarray}}
\newcommand{\ee}{\end{eqnarray}}

\def\ddbar{\partial\bar\partial}

\definecolor{Cyan}{cmyk}{1.,0,0,0}
\definecolor{Magenta}{cmyk}{0,1.,0,0}
\definecolor{Yellow}{cmyk}{0,0,1.,0}
\definecolor{White}{cmyk}{0,0,0,0}
\definecolor{Orange}{cmyk}{0,0.61,0.87,0}
\definecolor{RedOrange}{cmyk}{0,0.77,0.87,0}
\definecolor{Red}{cmyk}{0,1.,1.,0}
\definecolor{Purple}{cmyk}{0.45,0.86,0,0}
\definecolor{Violet}{cmyk}{0.79,0.88,0,0}
\definecolor{Blue}{cmyk}{1,0.5,0,0}
\definecolor{ProcessBlue}{cmyk}{0.96,0,0,0}
\definecolor{GreenYellow}{cmyk}{0.6,0,1.,0}
\definecolor{Black}{cmyk}{0,0,0,1}

\def\ric{\mathrm{Ric}}

\newtheorem{theorem}{Theorem}
\newtheorem{lemma}{Lemma}

\begin{document}
\starttext
\setcounter{footnote}{0}

\begin{center}

{\Large \bf On the K\"ahler-Ricci flow on Fano manifolds
\footnote{Work supported in part by the National Science Foundation under grants DMS-1855947 and DMS-1945869.}}

\vskip 0.2in

{\large Bin Guo, Duong H. Phong, and Jacob Sturm}

\vskip 0.3in

\begin{abstract}

A short proof of the convergence of the K\"ahler-Ricci flow on Fano manifolds admitting a K\"ahler-Einstein metric or a K\"ahler-Ricci soliton is given, using a variety of recent techniques. 

\end{abstract}
\end{center}

\setcounter{tocdepth}{2} 

\baselineskip=15pt
\setcounter{equation}{0}
\setcounter{footnote}{0}

\section{Introduction} 
\setcounter{equation}{0}
\label{sec:1}

In 2002, G. Perelman announced that he could show the convergence of the K\"ahler-Ricci flow on Fano manifolds admitting a K\"ahler-Einstein metric. This led to many attempts to supply a detailed proof, including \cite{TZ2007, TZ2013, TZZZ2013} and \cite{CS2016}. Perhaps the most transparent proof is that of Collins and Sz\'ekelyhidi \cite{CS2016}, which incorporated ideas from \cite{TZ2013} and a general version allowing a twisting. The understanding of the K\"ahler-Ricci flow has progressed considerably since that time, and many new ideas and techniques have been introduced. The purpose of this brief note is to point out that, if one combines the ideas and techniques of Darvas-Rubinstein \cite{DR},
Boucksom-Eyssidieux-Guedj \cite{BEG}, 
Collins-Sz\'ekelyhidi \cite{CS2016}, Kolodziej \cite{Ko}, and Phong-Song-Sturm-Weinkove \cite{PSSW, PSSW1, PSSW2}, a relatively short proof can now be written down. Thus we shall prove the following:

\begin{theorem} Let $(M,g)$ be a compact K\"ahler manifold with $c_1(M)>0$.
If $g$ is K\"ahler-Einstein, then the
K\"ahler-Ricci flow with any initial metric in $c_1(M)$   will converge 
smoothly and exponentially fast to a K\"ahler-Einstein metric. More generally, if $g$ is a
K\"ahler-Ricci soliton with vector field $X$, then the K\"ahler-Ricci flow with any initial metric in $c_1(M)$ which is invariant under Im $X$  will converge 
smoothly and exponentially fast to a K\"ahler-Ricci soliton, after modified by the one-parameter group generated by $X$. 
\end{theorem}
\noindent
Remark: Compared to the works \cite{TZZZ2013, CS2016}, our theorem provides a more precise convergence property of the K\"ahler-Ricci flow when there is a K\"ahler-Ricci soliton, namely, the gauge transformation is simply given by the group generated by $X$. Finally we remark that the theorem still holds when the initial metric is not necessarily invariant under Im $X$, by the recent work of Dervan-Sz\'ekelyhidi \cite{DS}.

\vskip .1in

\section{Proof of the Theorem}
\setcounter{equation}{0}

It may be instructive to compare the proof below to that of \cite{TZ2007} which treats the simpler case of Fano manifolds admitting a K\"ahler-Einstein metric and no non-trivial holomorphic vector fields.  There one develops a relative version of Kolodziej's capacity theory \cite{Ko} to get a one-sided bound on the potential. Then, combining the Moser-Trudinger \cite{T1997, PSSW0} inequality with a Harnack inequality, one gets the full $C^0$ bound. The higher order estimates are obtained by applying parabolic versions of Yau's estimates \cite{Y}.

\medskip

The present proof also begins with the Moser-Trudinger inequality, which is now available for K\"ahler-Einstein or K\"ahler-Ricci soliton manifold with holomorphic vector fields \cite{DR}. But we obtain $C^\alpha$ estimates for the potentials  directly from recent arguments introduced in \cite{GPS} using \cite{BEG, Ko}. Higher order estimates and convergence can be obtained now more efficiently by combining techniques from \cite{ST, PSSW2}.

\medskip
We begin by introducing the notation.
Let $(M,\omega_0)$ be a compact K\"ahler manifold with $c_1(M)>0$ and $X$ be a holomorphic vector field whose imaginary part $\mathrm{Im}\, X$ is a Killing vector field with respect to $\omega_0$. Write $\mathcal{K}_X$ for the space of K\"ahler metrics in $c_1(M)$ that are invariant under $\mathrm{Im}\, X$.  We write 
$$\mathcal P_X(M,\omega_0) = \{\varphi \in C^\infty(M,\mathbb R) ~|~ \omega_0 + i\ddbar \varphi \in \mathcal K_X, \, (\mathrm{Im}\, X) (\varphi) = 0\}$$ to be the space of K\"ahler potentials in $\mathcal K_X$. Given $\omega\in \mathcal{K}_X$, we define the Hamiltonian $\theta_{X,\omega}$ as the real-valued function satisfying 
$$\iota_X \omega = i \bar \partial \theta_{X,\omega},\quad \int_M e^{\theta_{X,\omega}}\omega^n = \int_M \omega^n = V.$$
We  define the Ricci potential $f = f_\omega$ by 
\begin{equation}\label{eqn:ric}-\ric(\omega) + \omega =i \ddbar f,\quad \int_M e^{-f}\omega^n = V.\end{equation}
We write $$u_{X,\omega} = f_\omega + \theta_{X,\omega}$$ to be the {\em modified Ricci potential}, which satisfies
\begin{equation}\label{eqn:ric potential}-\ric(\omega) + \omega + L_X \omega = i\ddbar u_{X,\omega},\end{equation}where $L_X$ denote the Lie derivative in the direction of $X$. It is known (\cite{TZ, PSSW2}) that if $\mathcal K_X\ni \omega  = \omega_0 + i\ddbar \varphi$ then $\theta_{X,\omega} = \theta_{X,\omega_0} + X\varphi$ and $\| \theta_{X,\omega}\|_{C^0}\le C(\omega_0, X)$ is uniformly bounded.

\smallskip

A K\"ahler metric $\omega_{KS}$ is a K\"ahler-Ricci soliton associated with $X$ if  $u_{X,\omega_{KS}} =0 $, i.e. 
$$\ric(\omega_{KS})  = \omega_{KS} + L_X \omega_{KS}.$$

The {\em modified Mabuchi $K$-energy} $\mu_{X,\omega_0}:\mathcal P_X(M,\omega_0)\to \mathbb R$ is defined by the variation (\cite{TZ})
$$\delta \mu_{X,\omega_0}(\varphi) = - \frac 1 V \int_M \delta\varphi (R-n - \nabla_j X^j - X(u_{X,\omega})  ) e^{\theta_{X,\omega}}\omega^n,\quad \mu_{X,\omega_0}(0) = 0$$
where $\omega = \omega_0+i\ddbar \varphi$ in the integral.

\smallskip

We will consider the following {\em normalized K\"ahler-Ricci flow} (\cite{PSSW2}, \cite{TZ})
\begin{equation}\label{eqn:mKRF}
\frac{\partial\omega}{\partial t} = - \ric(\omega) + \omega ,\quad \omega(0) = \omega_0,
\end{equation}
where $\omega_0\in\mathcal K_X$ is a fixed K\"ahler metric in $c_1(M)$. It is well-known that the solution $\omega(t)$ exists and lies in $ \mathcal K_X$ for all $t\in [0,\infty)$. 

We recall Perelman's estimates along the flow \eqref{eqn:mKRF}.

\begin{lemma}[\cite{ST}]\label{lemma 1}
There exists a constant $C=C(M, \omega_0)>0$ independent of $t$ such that
$$\|f_\omega\|_{C^0} + \| \nabla_\omega f_\omega\|_{C^0} + \|\Delta_\omega f_\omega\|_{C^0}  \le C$$ 
where $\omega = \omega(t)$ is the K\"ahler metric along the flow \eqref{eqn:mKRF}.

\end{lemma}

Let $Aut(M)$ be the automorphism group of $M$ and $\mathfrak{aut}(M)$ be the Lie algebra, i.e. the space of all real vector fields $X$ such that $\mathcal{L}_XJ=0$, where $J$ is the complex structure on $M$,  and denote $$\mathfrak{aut}_X(M) = \{Y\in \mathfrak{aut}(M)~|~ L_X Y = [X,Y] = 0\}.$$
Write $G=Aut_X(M)$ to be the connected Lie group whose Lie algebra is $\mathfrak{aut}_X(M)$. It follows that for any $\sigma\in G$, $d \sigma(X) = X$. 

For any $\sigma\in G$, we define $\sigma\cdot 0$ to be the K\"ahler potential of $\sigma^* \omega _0 = \omega_0 + i\ddbar \sigma\cdot 0$ (suitably normalized). Furthermore for any $\varphi\in \mathcal P_X(M,\omega_0)$, define $$\sigma\cdot\varphi = \varphi\circ \sigma + \sigma\cdot 0$$ which is the K\"ahler potential of $\sigma^* \omega$, where $\omega = \omega_0 + i\ddbar \varphi$.

\medskip

We recall the following Moser-Trudinger inequality from \cite{DR}. It is the extension to the case of K\"ahler-Ricci soliton manifolds with holomorphic vector fields of the sharp version for K\"ahler-Einstein manifolds with $\mathfrak{aut}(M)=0$ in \cite{PSSW0} of the inequality originally proved in \cite{T1997}:

\begin{lemma}[Theorem 8.1 in \cite{DR}]\label{lemma MT}
If $M$ admits a K\"ahler-Ricci soliton metric associated to the holomorphic vector field $X$, then there exist positive constants $C$ and $D$ depending only on $\omega_0$ and $X$ such that
$$\mu_{X,\omega_0}(\varphi) \ge C \mathcal J_G(G\varphi) - D,\quad \forall \varphi \in \mathcal P_X(M,\omega_0),$$
where $\mathcal J_G(G\varphi)= \inf_{\sigma \in G} J_{\omega_0}(\sigma\cdot\varphi)$ and $J_{\omega_0}$ is the Aubin-Yau $J$-functional with the reference metric $\omega_0$.
\end{lemma}

The following lemma, which is a version of estimates in \cite{GPS}, is the key starting point:

\begin{lemma}
\label{C-alpha}
Let $\{t_j\}$ be an arbitrary sequence of times, and denote $\omega_j = \omega(t_j) = \omega_0 + i\ddbar \varphi_j$. Then there exists a uniform constant $C>0$ and $\sigma_j\in G$ with
\bea
\|\psi_j\|_{C^\alpha(M,\omega_0)}\leq C
\eea
if we set $\tilde \omega_j = \sigma_j^* \omega_j\in \mathcal K_X$, and $\tilde\omega_j=\omega_0+i\ddbar\psi_j$, ${\rm sup}_M\psi_j=0$. In particular,
\begin{equation}\label{eqn:cMA 3}
C^{-1}\omega_0^n\le \tilde \omega_j^n\le C \omega_0^n.
\end{equation}
\end{lemma}

\bigskip
{\it Proof of Lemma \ref{C-alpha}}. 
From Lemma \ref{lemma MT}, for each $j$ there exists a $\sigma_j\in G$ such that $\mathcal J_G(G\varphi_j) \ge J_{\omega_0}(\sigma_j \cdot \varphi_j) - 1$ and $J_{\omega_0}(\sigma_j\cdot \varphi_j) \le C^{-1} \mu_{X,\omega_0}(\varphi_j) + C^{-1} D+ 1\le C(\omega_0, X)$, since the modified Mabuchi $K$-energy $\mu_{X,\omega_0}$ is non-increasing along the flow. 
It follows by a straightforward calculation that $f_{\tilde \omega_j} = \sigma_j^* f_{\omega_j}$ so $f_{\tilde \omega_j}$ is also uniformly bounded by Lemma \ref{lemma 1}. From the equation \eqref{eqn:ric}, the metric $\tilde\omega_j$  satisfies the complex Monge-Amp\`ere equation
\begin{equation}\label{eqn:cMA}
\tilde\omega_j^n=(\omega_0 + i\ddbar \psi_j)^n = e^{ -( \psi_j - \sup_M \psi_j  ) + f_{\tilde \omega_j} - f_{\omega_0} + c_j  }\omega_0^n,
\end{equation}
where $c_j$ is a normalizing constant so that both sides have the same integral over $M$. It follows from the normalizing condition $\frac 1 V \int_M e^{-f_{\tilde\omega_j}} \tilde \omega_j^n = 1$ and equation \eqref{eqn:cMA} that
$$V= \int_{M} e^{ - (\psi_j - \sup_M \psi_j) - f_{\omega_0} + c_j  }\omega_0^n\ge e^{c_j} \int_M e^{-f_{\omega_0} }\omega_0^n = e^{c_j} V\quad\Rightarrow \quad c_j\le 0.$$
So
\begin{equation}\label{eqn:cMA 1}
(\omega_0 + i\ddbar \psi_j)^n \le C e^{ -(\psi_j - \sup_M \psi_j)   }\omega_0^n=C e^{ -\tilde\psi_j    }\omega_0^n
\end{equation}
where $C=C(\omega_0, X)>0$ depends on the bound on $f_{\tilde \omega_j}$ and $f_{\omega_0}$. $\tilde \psi_j = \psi_j - \sup_M \psi_j$ belongs to a set $S_A$ for some $A = A(\omega_0, X)>0$, which is
$$S_A : = \{\psi\in \mathcal P_X (M,\omega_0)\subset \mathcal E_X^{1}(M,\omega_0)~|~ \sup_M \psi = 0\text{ and } J_{\omega_0} (\psi)\le A\}.$$
$S_A$ is compact under the weak $L^1(M,\omega_0^n)$-topology in $\mathcal E_X^1(M,\omega_0)$, and each $\psi\in S_A$ has zero Lelong numbers at any point $x\in M$. This implies, as in \cite{BEG}, that for any $p>0$, there exists a constant $C_p = C(\omega_0, A, p)>0$ such that 
\begin{equation}\label{eqn:cMA 2}\int_M e^{-p \psi}\omega_0^n \le C_p,\quad \forall \psi\in S_A.\end{equation}
Combining equations \eqref{eqn:cMA 1} and \eqref{eqn:cMA 2} and Kolodziej's theorem (\cite{Ko}), for $p>1$, there exist an $\alpha = \alpha(n, p)\in (0,1)$ and $C=C(\omega_0, A, p)>0$ such that 
\begin{equation}\label{eqn:estimate psi}
\| \tilde \psi\|_{C^\alpha(M,\omega_0)}\le C.
\end{equation}
The estimate (\ref{eqn:cMA 3}) follows next from the equation (\ref{eqn:cMA}). The lemma is proved.

\bigskip
We can give now the proof of the theorem. 
Let $\tilde \omega_j(s): = \sigma_j^*\omega(t_j + s)$ for $s\in [0,3]$. Then $\tilde \omega_j(s)$ satisfies the {\em K\"ahler-Ricci flow} equation
\begin{equation}\label{eqn:mKRF 2}
\frac{\partial \tilde \omega_j(s)}{\partial s} = -\ric(\tilde\omega_j(s)) + \tilde \omega_j(s) ,\quad \tilde \omega_j(0) = \tilde \omega_j.
\end{equation}
Let $\tilde \omega_j(s) = \omega_0 + i \ddbar \tilde \psi_j(s)$. The equation \eqref{eqn:mKRF 2} is equivalent to the following complex Monge-Amp\`ere equation for $\tilde \psi_j(s)$
\begin{equation}
\label{eqn:cMA 2}
\frac{\partial \tilde \psi_j(s)}{\partial s} = \log\frac{ (\omega_0 + i\ddbar \tilde\psi_j(s))^n  }{\omega_0^n} + \tilde \psi_j(s)  + f_{\omega_0},\quad \tilde \psi_j(0) = \tilde \psi_j,
\end{equation}
where $\tilde \psi_j$ is the K\"ahler potential of the initial metric $\tilde \omega_j$ with $\sup_M \tilde \psi_j = 0$. Equations \eqref{eqn:cMA 3} and \eqref{eqn:estimate psi} imply that 
$$\| \tilde \psi_j(0)\|_{C^0} + \| \frac{\partial}{\partial s}\Big|_{s=0} \tilde \psi_j(s)\|_{C^0}\le C(\omega_0,X).$$ 
We recall the following parabolic version of Yau's estimates for complex Monge-Amp\`ere equations.
\begin{lemma}[Proposition 2.1 in \cite{SzTo}]\label{lemma ST}
Suppose $\varphi\in PSH(M,\omega_0)$ satisfies the parabolic Monge-Amp\`ere equation $$\frac{\partial \varphi}{\partial t} = \log\frac{(\omega_0 + i\ddbar \varphi)^n}{\omega_0^n} + \varphi + f_0,\quad \varphi(0) = \varphi_0$$ for some smooth function $f_0: M\to \mathbb R$. Then for any $k\in\mathbb N$, there exists a smooth function  $C_k: (0, \infty) \to \mathbb R_+$ which depends only on $\| \varphi_0\|_{C^0}$ and $\| \dot\varphi|_{t=0}\|_{C^0}$ such that 
$$\| \varphi(t)\|_{C^k(M,\omega_0)}\le C_k(t),\quad \forall t\in (0,\infty).$$
\end{lemma}

We now apply Lemma \ref{lemma ST} to conclude that for each $k\in\mathbb N$, there exists a constant $C_k=C(\omega_0, X)>0$ such that 
$$\|\tilde  \psi_j(s)\|_{C^k(M,\omega_0)} \le C_k,\quad \forall s\in [1,2].$$
In particular, this implies the metric $\sigma_j^*(\omega(t_j+1))$ has a uniformly bounded K\"ahler potential and its derivatives are also uniformly bounded.

\medskip

 Since $t_j>0$ is arbitrarily chosen, we conclude that for any $t\ge 1$, there exists a $\sigma_t\in G$ such that 
$$\sigma_t^* ( \omega(t) ) = \omega_0 + i\ddbar \psi_t,\quad\text{with } \| \psi_t\|_{C^k(M,\omega_0)}\le C_k(\omega_0, X).$$

Let $\eta_t = \exp(tX)$ be the one-parameter group generated by $X$, and $\lambda(t)$ be the first nonzero eigenvalue of the $\bar \partial$-operator associated to $\eta_t^*(\omega(t))$. Since $\lambda(t)$ is also the corresponding eigenvalue of $\sigma_t^*(\omega(t))$, which satisfies uniform estimates and hence forms a compact set, we conclude that
\begin{equation*}
\inf_{t\in [0,\infty)} \lambda(t) >0.\end{equation*}

On the other hand, since $M$ admits a K\"ahler-Ricci soliton,
by \cite{TZ} we have \begin{equation*}
\inf_{\varphi\in\mathcal P_X(M,\omega_0)}\mu_{X,\omega_0}(\varphi)>-\infty.\end{equation*} Applying Theorem 3 in \cite{PSSW2} we  conclude that $\hat \omega_t = \eta_t^*(\omega(t))$ converge exponentially fast to a K\"ahler-Ricci soliton metric, observing that $\eta_t^*(\omega(t))$ satisfies the {\em modified K\"ahler-Ricci flow} defined in \cite{PSSW2}.

The proof of the theorem is complete.

\section{Additional remarks}
\setcounter{equation}{0}

\medskip

We also observe that another proof of the convergence of the K\"ahler-Ricci flow can also be obtained using the ideas in \cite{PSSW1, PSSW2} and the recent result in \cite{GPS}. For example, if $M$ is K\"ahler-Einstein, then the $K$-energy is bounded from below, and by \cite{GPS}, the lowest strictly positive eigenvalue of the Laplacian on vector fields is uniformly bounded away from $0$ along the K\"ahler-Ricci flow. By the main result of \cite{PSSW1}, the flow converges.


\bigskip

\noindent Bin Guo, Department of Mathematics and Computer Sciences, Rutgers University, Newark, NJ  07102

Email: bguo@rutgers.edu

\medskip

\noindent Duong H. Phong, Department of Mathematics, Columbia University, New York, NY 10027

Email: phong@math.columbia.edu

\medskip

\noindent Jacob Sturm,  Department of Mathematics and Computer Sciences, Rutgers University, Newark, NJ  07102

Email: sturm@andromeda.rutgers.edu


\begin{thebibliography}{99}

{\small

\bibitem{BEG}
Boucksom, S., Eyssidieux, P., and V. Guedj, ``Introduction to the K\"ahler-Ricci flow", Lecture Notes in Math., 2086, Springer, 2013.

\bibitem{CS2016}
Collins, T., and G. Sz\'ekelyhidi, 
``The twisted K\"ahler-Ricci flow"
J. Reine Angew. Math. 716 (2016), 179--205.
53C44 (32Q20 53C55)

\bibitem{DR} Darvas T. and Y. Rubinstein. ``Tian's properness conjectures and Finsler geometry of the space of K\"ahler metrics''. J. Amer. Math. Soc. 30 (2017), 347 -- 387.

\bibitem{DS} Dervan, R. and G. Sz\'ekelyhidi, ``The K\"ahler-Ricci flow and optimal degenerations'', arXiv: 1612.07299


\bibitem{GPS} Guo, B., D.H. Phong, and J. Sturm,
``K\"ahler-Einstein metrics and eigenvalue gaps", arXiv:2001.05794


\bibitem{Ko} Kolodziej, S. ``H\"older continuity of solutions to the complex Monge-Amp\`ere equation with the right hand side in $L^p$. The case of compact K\"ahler manifolds''. Math. Ann. 342, 379 -- 386 (2008)

\bibitem{N} Nadel, A.,
``Multiplier ideal sheaves and K\"ahler-Einstein metrics of positive scalar
curvature''. Ann. of Math. 132 (1990) 549 -- 596.

\bibitem{PSSW0} Phong, D. H., J. Song, J. Sturm, and B. Weinkove, ``The Moser-Trudinger inequality on K\"ahler-Einstein manifolds''
Amer. J. of Math. 130 no 4 (2008) 1067-1085.

\bibitem{PSSW} Phong, D.H., J. Song, J. Sturm, and B. Weinkove, ``The K\"ahler-Ricci flow with positive bisectional curvature''. Invent. Math. 173 (2008), no. 3, 651 - 665.

\bibitem{PSSW1} Phong, D. H., J. Song, J. Sturm, and B. Weinkove, ``The K\"ahler-Ricci flow and the $\bar\partial$ operator on vector fields''. J. Differential Geom. 81 (2009), no. 3, 631 -- 647.

\bibitem{PSSW2} Phong, D. H., J. Song, J. Sturm, and B. Weinkove, 
``The modified K\"ahler-Ricci flow and solitons", arXiv:0809.0941, 2008.

\bibitem{PS07} Phong, D.H. and J. Sturm, ``Lectures on stability and constant scalar curvature''. Handbook of geometric analysis, No. 3, 357 -- 436, Adv. Lect. Math. (ALM), 14, Int. Press, Somerville, MA, 2010.

\bibitem{ST} Sesum, N. and G. Tian, ``Bounding scalar curvature and diameter along the K\"ahler Ricci flow (after Perelman)''. J. Inst. Math. Jussieu 7 (2008), no. 3, 575 -- 587

\bibitem{Siu} Siu, Y.T., ``The existence of K\"ahler-Einstein metrics on manifolds with positive anticanonical line bundle and a suitable finite symmetry group''. Ann. of Math. 127 (1988) 585 -- 627.

\bibitem{SzTo} Sz\'ekelyhidi, G. and V. Tosatti, ``Regularity of weak solutions of a complex Monge-Amp\`ere equation'', Anal. PDE 4 (2011), no. 3, 369 -- 378

\bibitem{T1997} Tian, G., ``K\"ahler-Einstein metrics with positive scalar curvature''. Inventiones Math. 130 (1997) 1 -- 37.


\bibitem{TZ}Tian G. and X. H. Zhu, ``A new holomorphic invariant and uniqueness of K\"ahler-Ricci solitons''. Comment. Math. Helv. 77 (2002), 297 -- 325.

\bibitem{TZ2007} Tian G. and X. H. Zhu,
``Convergence of K\"ahler-Ricci flow". J. Amer. Math. Soc. 20 (2007), no. 3, 675--699.

\bibitem{TZ2013} Tian G. and X. H. Zhu,
``Convergence of the K\"ahler-Ricci flow on Fano manifolds". J. Reine Angew. Math. 678 (2013), 223--245.

\bibitem{TZZZ2013}Tian G., Zhang, S., Zhang, Z., and X. H. Zhu,
``Perelman's entropy and K\"ahler-Ricci flow on a Fano manifold". Trans. Amer. Math. Soc. 365 (2013), no. 12, 6669--6695.

\bibitem{Y} Yau, S.T., ``On the Ricci curvature of a compact K\"ahler manifold and the complex
Monge-Amp\`ere equation. I."
Comm. Pure Appl. Math. 31 (1978), no. 3, 339--411

}


\end{thebibliography}
\end{document}